\documentclass[10pt]{amsart}
\usepackage{amscd}
\usepackage{amssymb}

\newtheorem{thm}[equation]{Theorem}
\newtheorem{cor}[equation]{Corollary}
\newtheorem{prop}[equation]{Proposition}
\newtheorem{lem}[equation]{Lemma}
\theoremstyle{definition}
\newtheorem{dfn}[equation]{Definition}
\newtheorem{rem}[equation]{Remark}
\newtheorem{exa}[equation]{Example}

\numberwithin{equation}{section}

\newcommand{\iso}{\xrightarrow{\simeq}}
\newcommand{\inj}{\hookrightarrow}
\newcommand{\surj}{\twoheadrightarrow}
\newcommand{\xar}{\xrightarrow}
\newcommand{\opn}{\operatorname}
\newcommand{\opnt}[1]{\mathrm{#1}} 
\newcommand{\cat}[1]{\operatorname{\mathsf{#1}}}

\newcommand{\rmitem}[1]{\item[\text{\textup{(#1)}}]}
\newcommand{\mfrak}[1]{\mathfrak{#1}}

\newcommand{\mrm}[1]{\mathrm{#1}}
\newcommand{\mbb}[1]{\mathbb{#1}}

\newcommand{\tup}[1]{\textup{#1}}
\newcommand{\bsym}[1]{\boldsymbol{#1}}

\newcommand{\bra}[1]{\langle #1 \rangle}
\DeclareMathSymbol{\mathbbk}{\mathord}{AMSb}{"7C}
\renewcommand{\k}{\mathbbk}

\newcommand{\abs}[1]{\lvert #1 \rvert}

\author{Amnon Yekutieli}
\title{On the Structure of Behaviors}
\address{Department of  Mathematics 
Ben Gurion University, 
Be'er Sheva 84105, 
Israel}
\email{amyekut@math.bgu.ac.il}
\date{28 May 2004}
\thanks{
Primary: 93C05; Secondary: 93B25, 16W80, 13J10, 16P40.}

\begin{document}

\begin{abstract}
A behavior is a closed shift invariant subspace of the space of 
sequences with entries in a field $\k$. We work out an explicit 
duality for $\k$-modules. This duality is then used to derive 
properties of behaviors, and their noncommutative generalizations.
\end{abstract}

\maketitle

\setcounter{section}{-1}
\section{Introduction}

Let $\k$ be a field. Put on $\k$ the discrete topology, and put on 
\[ \k^{\mbb{N}} = \{ \phi : \mbb{N} \to \k \} = 
\prod_{i \in \mbb{N}} \, \k \]
the product topology. 
Then $\k^{\mbb{N}}$ is a topological $\k$-module. 
The shift operator $\sigma : \k^{\mbb{N}} \to \k^{\mbb{N}}$ is
\[ \sigma(\lambda_0, \lambda_1, \lambda_2, \ldots) :=
(\lambda_1, \lambda_2, \ldots) . \]
A {\em behavior} (or {\em discrete linear system}) 
is a closed shift invariant $\k$-submodule
$M \subset (\k^{\mbb{N}})^r$ for some natural number $r$. 
To be precise, this is a {\em $1$-dimensional behavior}. 
For any $n \geq 1$ there is a corresponding definition of 
{\em $n$-dimensional behavior},
which is a closed shift invariant $\k$-submodule of
$(\k^{\mbb{N}^n})^r$. But in the introduction we will stick to the 
$1$-dimensional case. The notion of behavior is attributed to 
J.C. Willems. See \cite{Fu1} 
and its references for more background information. 

A behavior is naturally a module over the polynomial ring $\k[z]$; 
the variable $z$ acts by the shift operator $\sigma$. 
Early on it was realized that behaviors are intimately related to 
finitely generated $\k[z]$-modules, and
that in fact these are dual mathematical objects. Yet some of the 
more subtle aspects were not fully understood (cf.\ \cite{Fu2}). 

There does exist a detailed treatment of behaviors by Oberst
in \cite{Ob}. However we feel that this treatment is unduly 
complicated; perhaps because the author wanted to consider 
discrete and continuous systems in a unified fashion. Thus Oberst 
showed that the $\k[z]$-module $\k^{\mbb{N}}$ is an 
injective cogenerator of the category $\cat{Mod} \k[z]$ of 
$\k[z]$-modules, and then he considered behaviors as 
modules over the endomorphism 
ring $\opn{End}_{\k[z]}(\k^{\mbb{N}})$. A similar approach was 
taken in \cite{Lo}. In doing so some interesting features of the 
theory were missed (such as the counterexample in Section 5, or
the noncommutative generalization in Section 6).

The aim of our paper is to clarify the structure of behaviors and 
to point out possible generalizations. 

Our basic observation is that the duality underlying behaviors has 
nothing to do with polynomial rings -- it is a duality for 
$\k$-modules (i.e.\ vector spaces). We establish a duality 
$\mrm{D} : M \mapsto M^*$ between 
the category $\cat{Mod} \k$ of $\k$-modules, and the category 
$\cat{TopMod}_{\mrm{pf}} \k$ of {\em profinite topological 
$\k$-modules}. By definition a profinite topological 
$\k$-module is a topological $\k$-module that is an inverse 
limit of finitely generated 
$\k$-modules (with discrete topologies). A morphism in 
$\cat{TopMod}_{\mrm{pf}} \k$ is a continuous $\k$-linear 
homomorphism. We prove that a closed $\k$-submodule of a profinite 
topological $\k$-module is also profinite.
It turns out (as a consequence of the duality)
that the profinite topological $\k$-modules are 
precisely the linearly compact topological vector spaces, in the 
sense of \cite[Section 10.9]{Ko}.  

Of course this duality is just a very easy instance of 
Gabriel-Matlis duality, valid for noncommutative noetherian 
rings and suitable module categories over them. But instead 
of quoting from classical (and complicated) work in ring theory, 
or alternatively dressing up the more naive duality results 
of \cite{Ko}, we simply work out the proofs explicitly. 
This provides us with a 
lot of information that is particular to $\k$-modules. 

Given a commutative $\k$-algebra $A$, we consider the category
$\cat{TopMod}_{\mrm{pf} / \k} A$ 
consisting of topological $A$-modules that are profinite as 
$\k$-modules. Then the duality $\mrm{D}$ restricts to a duality
$\mrm{D} : \cat{Mod} A \to \cat{TopMod}_{\mrm{pf} / \k} A$.

In this framework the behaviors can be described as follows. Take 
$A := \k[z]$. The shift module $\k^{\mbb{N}}$ coincides with 
the dual $\k[z]^*$ as topological $\k[z]$-modules, and it is
profinite over $\k$. A behavior $M$ is then a closed 
$\k[z]$-submodule of  $(\k[z]^*)^r$ 
for some natural number $r$; and so it 
is profinite as $\k$-module. 
By duality $M^*$ is a quotient of $\k[z]^r$, and hence it is a
finitely generated $\k[z]$-module. This argument can be reversed. 
See Theorems \ref{thm2} and \ref{thm4} for precise statements
(for $n$-dimensional behaviors). 

The paper is organized as follows. In Section 1 we study profinite 
topological $\k$-modules. Section 2 is devoted to the duality 
between the category $\cat{Mod} \k$ of $\k$-modules and the 
category $\cat{TopMod}_{\mrm{pf}} \k$ of profinite 
topological $\k$-modules. In Section 3 we consider the category 
$\cat{TopMod}_{\mrm{pf} / \k} A$ for any $\k$-algebra (not 
necessarily commutative). In Section 4 we concentrate on 
$n$-dimensional behaviors, which are modules over 
the polynomial ring $\k[z_1, \ldots, z_n]$. 
In Section 5 we exhibit 
$\k[z]$-linear homomorphisms between behaviors that 
are not continuous. Finally in Section 6 we study noncommutative 
generalizations of behaviors.

\bigskip \noindent
\textbf{Acknowledgements.}
I wish to thank Paul Fuhrmann for getting me interested in
behaviors. Also thanks to Sverre Smal{\o} for reading the paper and 
suggesting some improvements. Finally thanks to the referees for 
their remarks.

\section{Profinite Topological $\k$-Modules}

In this section $\k$ is a noetherian commutative ring (e.g.\ a 
field or the ring of integers $\mbb{Z}$).
Denote by $\cat{Mod} \k$ the category of $\k$-modules. 

By a topological $\k$-module 
we mean a $\k$-module $V$ endowed with a topology 
(of any sort) such that addition $V \times V \to V$ is continuous, 
and for any $\lambda \in \k$ the multiplication 
map $\lambda : V \to V$ is continuous. Let $\cat{TopMod} \k$ be the 
category whose objects are the topological $\k$-modules and whose 
morphisms are the continuous $\k$-linear homomorphisms. So
\[ \opn{Hom}_{\cat{TopMod} \k}(W_1, W_2) = 
\opn{Hom}^{\mrm{cont}}_{\k}(W_1, W_2) . \]
The category $\cat{TopMod} \k$ is additive, but it is not abelian. 
See \cite[Chapter II]{HS} for background material on categories
and functors. 

Let $\phi : W_1 \to W_2$ be a morphism in $\cat{TopMod} \k$. The 
morphism $\phi$ is called a {\em strict monomorphism} if $\phi$ is 
injective and $W_1$ has the subspace topology induced from $W_2$. 
The morphism $\phi$ is called a
{\em strict epimorphism} if $\phi$ is surjective and $W_2$ has 
the quotient topology induced from $W_1$. And
$\phi$ is called a {\em strict morphism} if it factors into
$\phi = \phi_2 \circ \phi_1$,
where $\phi_1 : W_1 \to W$ is a strict epimorphism
and $\phi_2 : W \to W_2$ is a strict monomorphism.
A sequence of homomorphisms (possibly infinite on either side)
\[ \cdots \to W_0 \xar{\phi_0} W_1 \xar{\phi_1} W_2 \to \cdots \]
in $\cat{TopMod} \k$ is called {\em strict-exact} if for all $i$ 
one has $\opn{Im}(\phi_{i - 1}) = \opn{Ker}(\phi_i)$, and $\phi_i$ is
strict. See \cite{Bo} for more details on strict homomorphisms. 

Recall that a {\em quasi-ordered} set is a pair $I = (I_0, I_1)$ 
consisting of a set $I_0$, and a set of arrows 
$I_1  = \{ i \to j \} \subset I_0 \times I_0$, such that $I_1$
is a reflexive and transitive relation on $I_0$. Thus $I$ is a 
category with at most one morphism between any two objects. 
We note that the pair
$I^{\mrm{op}} := (I_0, I_1^{\mrm{op}})$, in which all the arrows 
are reversed, is also a quasi-ordered set. 
A quasi-ordered set $I$ is called a {\em directed set} if for
any two elements $i, j \in I_0$ there exist arrows
$i \to k$ and $j \to k$ in $I_1$. 

Suppose $I$ is a quasi-ordered set and $\cat{C}$ is any category.
A {\em direct system} in $\cat{C}$ indexed by $I$ is a functor 
$F : I \to \cat{C}$. Usually one refers to this by saying that
$\{ C_i \}_{i \in I}$ is a direct system, where for any
$i \in I_0$ one writes $C_i := F i \in \cat{C}$, and the morphisms
$F(i \to j) : C_i \to C_j$ are implicit. A {\em direct limit} of 
the system $\{ C_i \}_{i \in I}$ is an object $C \in \cat{C}$,
equipped with a compatible system of morphisms 
$\phi_i : C_i \to C$, which is universal for this property.
By `compatible system of morphisms' we mean that for every arrow 
$\alpha : i \to j$ in $I$ one has
$\phi_j \circ \alpha = \phi_i$. And by `universal' we mean that 
given any $D \in \cat{C}$ with a compatible system of morphisms 
$\psi_i : C_i \to D$ there is exactly one morphism 
$\psi : C \to D$ such that 
$\psi_i = \psi \circ \phi_i$. 
A direct limit is unique if it exists; and in many cases it does 
exist, e.g.\ when $\cat{C}$ is the category of sets $\cat{Sets}$ or 
the category $\cat{Mod} \k$. 
The direct limit $C$ is denoted by
$\varinjlim F$ or 
$\underset{i \in I}{\varinjlim}\, C_i$. 

An {\em inverse system} in $\cat{C}$ indexed by $I$ is a functor 
$F : I^{\mrm{op}} \to \cat{C}$. So there is a set 
$\{ C_i \}_{i \in I}$ of objects of $\cat{C}$, and for each arrow
$\alpha : i \to j$ in $I$ we get an arrow 
$F(\alpha) : C_j \to C_i$ in $\cat{C}$. An {\em inverse limit} 
is an object $C \in \cat{C}$, with a compatible system of 
morphisms $C \to C_i$, which is universal for this property. It is 
denoted by $\varprojlim F$ or 
$\underset{i \in I}{\varprojlim}\, C_i$. 

Observe that $F : I^{\mrm{op}} \to \cat{C}$ can be viewed as a 
functor
$F : I \to \cat{C}^{\mrm{op}}$, where $\cat{C}^{\mrm{op}}$ is the 
opposite category (same objects but reversed morphisms). So an 
inverse system in $\cat{C}$ is the same as a direct system in 
$\cat{C}^{\mrm{op}}$, and $\varprojlim F$ in $\cat{C}$ coincides 
with $\varinjlim F$ in $\cat{C}^{\mrm{op}}$.

Refer to \cite[Chapter 2]{Ro} for more on limits in categories.

We shall be mainly interested in inverse limits in 
$\cat{TopMod} \k$. Given an inverse system $\{ W_i \}_{i \in I}$ 
of topological $\k$-modules indexed by a quasi-ordered set $I$, 
the inverse limit $W := \underset{i \in I}{\varprojlim}\, W_i$
is constructed like this. Inside the 
product $\prod_{i \in I} W_i$, endowed with the product topology,
one takes the submodule consisting 
of sequences $(w_i)$ such that for any arrow $i \to j$ in $I$ the 
homomorphism $W_j \to W_i$ sends $w_j \mapsto w_i$. The limit module 
$W$ is equipped with a system of continuous $\k$-linear 
homomorphisms 
$\pi_i : W \to W_i$. Given any topological $\k$-module $U$ the 
system of morphisms $\{ \pi_i \}$ induces a bijection of sets
\[ \opn{Hom}^{\mrm{cont}}_{\k}(U, W) \to
\underset{i \in I}{\varprojlim}\, 
\opn{Hom}^{\mrm{cont}}_{\k}(U, W_i) , \]
where this last limit is taken in $\cat{Mod} \k$.
Cf.\ \cite[Section 1.1]{Ye1}. 

\begin{dfn} \label{dfn8}
Given a topological $\k$-module $W$ let us denote by
$\cat{Cofin} W$ the set of open cofinite $\k$-submodules of $W$, 
namely those submodules $W' \subset W$ such that the quotient 
$W / W'$ is a discrete finitely generated $\k$-module. 
\end{dfn}

The set $\cat{Cofin} W$ is partially ordered by inclusion, and 
hence it is a quasi-ordered set. 

\begin{lem} \label{lem7}
The quasi-ordered sets $\cat{Cofin} W$ and 
$(\cat{Cofin} W)^{\mrm{op}}$ are both directed. 
\end{lem}

\begin{proof}
First consider $\cat{Cofin} W$. Let $W_1, W_2 \in \cat{Cofin} W$. 
We have to show that there is some $W_3 \in \cat{Cofin} W$ such that 
$W_1 \subset W_3$ and $W_2 \subset W_3$. 
For this we can take $W_3 := W$.

Next consider $(\cat{Cofin} W)^{\mrm{op}}$. 
We have to show that there is some $W_3 \in \cat{Cofin} W$ such that 
$W_3 \subset W_1$ and $W_3 \subset W_2$. 
Take the intersection $W_3 := W_1 \cap W_2$. 
Clearly $W_3$ is open in $W$. It remains
to show it is cofinite. We know that $W / W_1$ and $W / W_2$ 
are finitely generated. Because $\k$ is noetherian the submodule
$\frac{W_1}{W_1 \cap W_2} \subset \frac{W}{W_2}$
is finitely generated. By the exact sequence
\[ 0 \to \frac{W_1}{W_1 \cap W_2} \to 
\frac{W}{W_1 \cap W_2} \to 
\frac{W}{W_1} \to 0 \]
it follows that $\frac{W}{W_1 \cap W_2}$ is finitely generated.
\end{proof}

\begin{dfn} \label{dfn0}
Let $W$ be a topological $\k$-module. For any 
$W' \in \cat{Cofin} W$ the quotient $W / W'$ has the discrete  
topology. The inverse limit $\varprojlim W / W'$, as $W'$ runs 
over the quasi-ordered set $\cat{Cofin} W$, 
is endowed with the $\varprojlim$ topology. 
The topological $\k$-module $W$ is called {\em profinite} if the 
canonical homomorphism
\[ W \to \underset{W' \in \cat{Cofin} W}{\varprojlim}\, W / W'  \]
is an isomorphism in $\cat{TopMod} \k$.
\end{dfn}

We denote by $\cat{TopMod}_{\mrm{pf}} \k$ the full subcategory of 
$\cat{TopMod} \k$ consisting of profinite modules. 

\begin{lem} \label{lem8}
Let $W$ be a profinite topological $\k$-module. Then the set
$\cat{Cofin} W$ is a basis of the topology at $0$; namely any open 
neighborhood $U$ of $0$ contains some open cofinite submodule 
$W' \subset W$. 
\end{lem}

\begin{proof}
For convenience let us write $I := \cat{Cofin} W$ and
$\cat{Cofin} W = \{ W'_i \}_{i \in I}$. 
By the definition of the product topology $U$ contains 
$W \cap \prod_{i \in I} U_i$ where $U_i \subset W / W'_i$
is open and $U_i = W / W'_i$ for all but finitely many $i$. Since
$W / W'_i$ is discrete we can assume that
$U_{i} = 0$ for all $i$ in some finite subset $I' \subset I$, 
and $U_i = W / W'_i$ for all $i \notin I'$. Now
\[ W \cap \prod_{i \in I} U_i = \bigcap_{i \in I'} 
\opn{Ker}(W \to W / W'_{i}) = 
\bigcap_{i \in I'}  W'_{i} . \]
But by Lemma \ref{lem7} the submodule
$W' := \bigcap_{i \in I'} W'_{i} \subset W$ is open and cofinite. 
\end{proof}

\begin{lem} \label{lem4}
Suppose $W$ is a profinite topological $\k$-module and $V$ is a 
discrete topological $\k$-module. Then the canonical map of sets
\[ \Psi : \underset{W' \in \cat{Cofin} W}{\varinjlim}\,
\opn{Hom}_{\k}(W / W', V) \to
\opn{Hom}^{\mrm{cont}}_{\k}(W, V)  \]
is bijective. In words, any continuous $\k$-linear homomorphism
$\phi : W \to V$ factors via some 
discrete finitely generated quotient $W / W'$. 
\end{lem}

\begin{proof}
Clearly $\Psi$ is injective. Suppose we are given a continuous 
$\k$-linear homomorphism $\phi : W \to V$. 
Because $\{ 0 \} \subset V$ is open the kernel 
$\opn{Ker}(\phi)$ is open in $W$. According to Lemma \ref{lem8}
there is some $W' \in \cat{Cofin} W$ such that 
$W' \subset \opn{Ker}(\phi)$. Then denoting by 
$\bar{\phi} : W / W' \to V$ the induced homomorphism, we see that 
$\phi = \Psi(\bar{\phi})$.
\end{proof}

\begin{lem} \label{lem6}
Let $W$ be a profinite topological $\k$-module. Then
$W$ is separated \tup{(}i.e.\ Hausdorff\tup{)}.
\end{lem}

\begin{proof}
This is because the product 
$\underset{W' \in \cat{Cofin} W}{\prod}\, W / W'$ 
is separated, and $W$ has the subspace topology.
\end{proof}

\begin{rem}
Actually a profinite topological $\k$-module $W$ is also complete 
(in the sense of Cauchy filters); but we do not need this fact. 
\end{rem}

\begin{prop} \label{prop4}
Suppose $\{ W_i \}_{i \in I}$ is an inverse system of 
finitely generated 
discrete $\k$-modules indexed by a quasi-ordered set $I$, 
and let
$W := \underset{i \in I}{\varprojlim}\, W_i$
in $\cat{TopMod} \k$. Assume $I^{\mrm{op}}$ is directed. 
Then $W$ is a profinite topological $\k$-module.
\end{prop}

\begin{proof}
Replacing $W_i$ by the discrete finitely generated module
$\opn{Im}(W \to W_i)$ we can assume the inverse 
system has surjections $W \surj W_i$ for all $i$. For each $i$ let
$W'_i := \opn{Ker}(W \to W_i)$, which is an open cofinite 
submodule. It suffices to prove that the inverse system
$\{ W'_i \}_{i \in I}$ is cofinal in the quasi-ordered set
$(\cat{Cofin} W)^{\mrm{op}}$. So pick any
$W' \in \cat{Cofin} W$. Since $W'$ is open we have
$\bigcap_{i \in I'} W_{i}' \subset W'$ for some finite subset
$I' \subset I$. Now $I^{\mrm{op}}$ is directed, so 
there is some $j \in I$ with arrows $j \to i$ for all $i \in I'$. 
Hence $W_j' \subset \bigcap_{i \in I'} W_{i}'$. 
\end{proof}

\begin{thm} \label{thm0}
Let $\k$ be a noetherian commutative ring, let
$W$ be a profinite topological $\k$-module and let $U$ be a 
closed submodule. Then $U$, with the induced topology, is also a 
profinite topological $\k$-module.
\end{thm}

\begin{proof}
As before let $I$ denote the quasi-ordered set
$\cat{Cofin} W$ and write $W'_i$ for the submodule labelled 
$i \in I$. Define $W_i := W / W'_i$ with the discrete topology.
By definition $W \cong \underset{i \in I}{\varprojlim}\, W_i$
in $\cat{TopMod} \k$. Let $U_i := \opn{Im}(U \to W_i)$, with the 
discrete topology. Since $\k$ is noetherian the module $U_i$ is 
finitely generated. By Proposition \ref{prop4} the inverse limit 
$\bar{U} := \underset{i \in I}{\varprojlim}\, U_i$ is a 
profinite topological $\k$-module. Since $U_i \subset W_i$ for all 
$i$ we get an injection $\bar{U} \inj W$. By 
Lemma \ref{lem8} and because $U \to U_i$ is surjective for every $i$
it follows that $U$ is dense in 
$\bar{U}$. But $U$ is closed in $W$, so we get $U = \bar{U}$. 
\end{proof}

\begin{rem}
If $W \cong \underset{i \in \mbb{N}}{\varprojlim}\, W_i$ 
for some inverse system $\{ W_i \}_{i \in \mbb{N}}$ of discrete 
modules, and $U \subset W$ is a closed submodule, 
then the quotient $W / U$ 
is also profinite; see \cite[Proposition 1.1.6]{Ye1}.
(By the way in this case $W$ is a metrizable topological space.) 
We do not know if $W / U$ 
is profinite when $W$ is not a countable inverse limit
of discrete modules.
\end{rem}

\section{Duality for Topological $\k$-Modules}

 From here on $\k$ is a field. 

Given a $\k$-module $V$ let 
$\cat{Fin} V$ be the set of finitely generated 
$\k$-submodules of $V$. 
Inclusion makes it a partially ordered set. It is directed, since 
given $V_1, V_2 \in \cat{Fin} V$ the sum $V_1 + V_2$ is also in 
$\cat{Fin} V$. And $V$ is the direct limit 
\[ V = \underset{V' \in \cat{Fin} V}{\varinjlim} V' \, . \]

\begin{lem} \label{lem5}
Given a $\k$-module $V$ the inclusions $V' \to V$ induce a 
bijection
\[ \Psi : \opn{Hom}_{\k}(V, \k) \iso
\underset{V' \in \cat{Fin} V}{\varprojlim} \opn{Hom}_{\k}(V', \k) 
. \]
\end{lem}

\begin{proof}
Say $\phi : V \to \k$ is nonzero. Then $\phi|_{V'} \neq 0$ for 
some $V' \in \cat{Fin} V$. This shows $\Psi$ is injective. 
Conversely, suppose $\{ \phi|_{V'} : V' \to \k \}$ is a compatible 
system of homomorphisms. Then we can patch the homomorphisms 
$\phi_{V'}$ to a ``global'' homomorphism $\phi : V \to \k$. So 
$\Psi$ is surjective. 
\end{proof}

\begin{dfn} \label{dfn6}
Given a $\k$-module $V$ let 
$\opnt{D} V := \opn{Hom}_{\k}(V, \k)$ 
be its dual module. We shall make 
$\opnt{D} V$ into a topological $\k$-module as follows.
\begin{enumerate}
\item If $V$ is finitely generated then $\opnt{D} V$ has the 
discrete topology.
\item For any $V$ we put on 
\[ \opnt{D} V \cong 
\underset{V' \in \cat{Fin} V}{\varprojlim} \opnt{D} V' \]
the inverse limit topology.
\end{enumerate}
\end{dfn}

\begin{lem} \label{lem1}
The assignment $V \mapsto \opnt{D} V$ gives rise to a contravariant 
functor
\[ \opnt{D} : \cat{Mod} \k \to \cat{TopMod} \k . \]
\end{lem}

\begin{proof}
We must show that for any $\k$-linear homomorphism
$\phi : V_1 \to V_2$ there is an induced continuous homomorphism
$\opnt{D}(\phi) : \opnt{D} V_2 \to \opnt{D} V_1$. 
Now let $V_1' \in \cat{Fin} V_1$. Then 
$V_2' := \phi(V_1') \in \cat{Fin} V_2$, and we get a continuous 
homomorphism
\[ \opnt{D}(\phi|_{V_1'}) : \opnt{D} V_2 \to \opnt{D} V_2'
\to \opnt{D} V_1' . \] 
Passing to the inverse limit in $V_1'$ we get continuous 
homomorphism
$\opnt{D}(\phi) : \opnt{D} V_2 \to \opnt{D} V_1$.
\end{proof}

\begin{prop} \label{prop1}
Let $V \in \cat{Mod} \k$. 
\begin{enumerate}
\item The topological $\k$-module $\opnt{D} V$ is profinite.
\item  The topology on $\k$-module $\opnt{D} V$ coincides with the 
weak$^*$ topology. 
\end{enumerate}
\end{prop}

\begin{proof}
(1) See Proposition \ref{prop4}.

\medskip \noindent 
(2) Recall that in the weak$^*$ topology on 
$W := \opn{Hom}_{\k}(V, \k) = \opnt{D} V$ 
a fundamental system of neighborhoods 
of $0$ is the set
\[ \{ W(v_1, \ldots, v_n) \mid n \geq 0,\ 
v_1, \ldots, v_n \in V \} , \]
where
\[ W(v_1, \ldots, v_n) := \{ \phi \in W \mid 
\phi(v_1) = \cdots = \phi(v_n) = 0 \} . \]
On the other hand a fundamental system of neighborhoods 
of $0$ in the topology of Definition \ref{dfn6} consists of the 
set
\[ \{ \opn{Ker}(W \to \opnt{D} V') \mid V' \in \cat{Fin} W \} . \]
Now letting
$V' := \sum_{i = 1}^n \k v_i \subset V$ we have
\[  W(v_1, \ldots, v_n) = 
\opn{Ker}(W \to \opnt{D} V') . \]
Since any finitely generated submodule $V'$ arises in this way, 
the two topologies coincide. 
\end{proof}

Take an object $W \in \cat{TopMod} \k$. Denote by 
\[ \opnt{D}^{\mrm{c}} W := \opn{Hom}_{\k}^{\mrm{cont}}(W, \k) =
\opn{Hom}_{\cat{TopMod} \k}(W, \k) \]
the continuous dual, which we consider as a $\k$-module 
(without any topology). Thus we get a contravariant functor
\[ \opnt{D}^{\mrm{c}} : \cat{TopMod} \k \to \cat{Mod} \k . \]

Given $V \in \cat{Mod} \k$ there is an adjunction (or evaluation) 
homomorphism 
$\alpha_V : V \to \opnt{D}^{\mrm{c}} \opnt{D} V$ 
whose formula is $\alpha_V(v)(\phi) := \phi(v)$. 
Likewise given $W \in \cat{TopMod} \k$ there is an adjunction 
homomorphism 
$\beta_W : W \to \opnt{D} \opnt{D}^{\mrm{c}} W$ 
whose formula is \linebreak
$\beta_W(w)(\psi) := \psi(w)$. 

Observe that a finitely generated $\k$-module $V$ is an object of 
$\cat{Mod} \k$, and at the same time, when endowed with the 
discrete topology, $V$ is an object of 
$\cat{TopMod}_{\mrm{pf}} \k$. Such $V$ is clearly reflexive, in 
the sense that $\alpha_V : V \to \opnt{D}^{\mrm{c}} \opnt{D} V$ 
is an isomorphism. 

It turns out that {\em all $\k$-modules are reflexive}:

\begin{thm} \label{thm5}
\begin{enumerate}
\item Let $V \in \cat{Mod} \k$. Then the adjunction homomorphism
$\alpha_V : V \to \opnt{D}^{\mrm{c}} \opnt{D} V$ 
is an isomorphism in $\cat{Mod} \k$.
\item Let $W \in \cat{TopMod}_{\mrm{pf}} \k$. 
Then the adjunction homomorphism
$\beta_W : W \to \opnt{D} \opnt{D}^{\mrm{c}} W$ 
is an isomorphism in $\cat{TopMod}_{\mrm{pf}} \k$.
\item The functor 
$\opnt{D} : \cat{Mod} \k \to \cat{TopMod}_{\mrm{pf}} \k$
is a duality \tup{(}i.e.\ a contravariant equivalence\tup{)}, 
with adjoint $\opnt{D}^{\mrm{c}}$.
\end{enumerate}
\end{thm}

\begin{proof}
(1) Define $I := \cat{Fin} V$, and rewrite
$\cat{Fin} V = \{ V'_i \}_{i \in I}$. Also let $W := \opnt{D} V$.
So $V \cong \underset{i \in I}{\varinjlim}\, V'_i$ and
$W = \underset{i \in I}{\varprojlim}\, \opnt{D} V'_i$. For any 
index $i$ let $W_i' := \opn{Ker}(W \to \opnt{D} V'_i)$, so 
$\opnt{D} V'_i \cong W / W_i'$.
The inverse system $\{ W_i' \}_{i \in I}$ is cofinal in 
$(\cat{Cofin} W)^{\mrm{op}}$, hence by Lemma \ref{lem4} we have
\[ \opnt{D}^{\mrm{c}} W \cong 
\underset{W' \in \cat{Cofin} W}{\varinjlim}\, 
\opnt{D}^{\mrm{c}} (W / W')
\cong \underset{i \in I}{\varinjlim}\, \opnt{D}^{\mrm{c}} (W / W'_i)
\cong \underset{i \in I}{\varinjlim}\, \opnt{D}^{\mrm{c}} 
\opnt{D} V'_i 
\cong \underset{i \in I}{\varinjlim}\, V'_i \cong V . \]
Denote by $\psi : \opnt{D}^{\mrm{c}} W \iso V$
the composition of this chain of isomorphisms, going from left 
to right. Then  
$\psi \circ \alpha_V = \bsym{1}_V : V \to V$, and therefore
$\alpha_V$ is an isomorphism.

\medskip \noindent
(2) Define $I := \cat{Cofin} W$, and rewrite
$\cat{Cofin} W = \{ W'_i \}_{i \in I}$. Let
$V := \opnt{D}^{\mrm{c}} W$, and for any index $i$ let
$V'_i := \opnt{D}^{\mrm{c}} (W / W'_i)$. 
 From Lemma \ref{lem4} we know that
$V \cong \underset{i \in I}{\varinjlim}\, V'_i$
in $\cat{Mod} \k$. It follows that the direct system
$\{ V'_i \}_{i \in I}$ is cofinal in $\cat{Fin} V$. 
Also $\opnt{D} V'_i \cong W / W'_i$. Hence
\[ \opnt{D} \opnt{D}^{\mrm{c}} W = \opnt{D} V =
\underset{V' \in \cat{Fin} V}{\varprojlim}\, \opnt{D} V' \cong
\underset{i \in I}{\varprojlim}\, \opnt{D} V'_i \cong
\underset{i \in I}{\varprojlim}\, W / W'_i \cong W \]
in $\cat{TopMod} \k$. 
Denote by $\phi : \opnt{D} \opnt{D}^{\mrm{c}} W \iso W$
the composition of this chain of isomorphisms, going from left 
to right. Then $\phi \circ \beta_W = \bsym{1}_W : W \to W$, and 
therefore $\beta_W$ is an isomorphism.

\medskip \noindent
(3) This follows from parts (1) and (2).
\end{proof}

In particular the theorem tells us that:
 
\begin{cor} \label{cor8}
For any $V_1, V_2 \in \cat{Mod} \k$ there is a canonical
$\k$-linear isomorphism
\[ \opn{Hom}_{\k}(V_1, V_1) \iso
\opn{Hom}_{\k}^{\mrm{cont}}(\opnt{D} V_2, \opnt{D} V_1) ,\
\phi \mapsto \opnt{D}(\phi) . \]
\end{cor}

For the sake of convenience, and when no confusion may arise,
we shall write $V^* := \opnt{D} V$ for $V \in \cat{Mod} \k$, 
and also $W^* := \opnt{D}^{\mrm{c}} W$
for $W \in \cat{TopMod} \k$. 

\begin{cor} \label{cor1}
Let $\phi : W_1 \to W_2$ be a morphism in 
$\cat{TopMod}_{\mrm{pf}} \k$. Then:
\begin{enumerate}
\item The morphism $\phi$ is strict.
\item The $\k$-modules $\opn{Ker}(\phi)$, $\opn{Im}(\phi)$ and
$\opn{Coker}(\phi)$, with their induced topologies, are profinite. 
\item The module $\opn{Im}(\phi)$ is closed in $W_2$. 
\end{enumerate}
\end{cor}

\begin{proof}
(1) Applying $\opn{D}^{\mrm{c}}$ we get 
$\phi^* : W_2^* \to W_1^*$ in $\cat{Mod} \k$. Define
$V_{2,2} := \opn{Ker}(\phi^*)$, and choose a complement
$W_2^* = V_{2,1} \oplus V_{2,2}$. Similarly let
$V_{1,1} := \opn{Im}(\phi^*)$, and choose a complement
$W_1^* = V_{1,1} \oplus V_{1,2}$. Thus $\phi^*$ has a matrix 
representation
\[ \phi^* : W_2^* = V_{2,1} \oplus V_{2,2} 
\xar{\cdot 
\left[ \begin{smallmatrix} \psi_{1,1} & 0 \\ 0 & 0 
\end{smallmatrix} \right]}
V_{1,1} \oplus V_{1,2} = W_1^* \]
where $\psi_{1,1}$ is an isomorphism. Dualizing back we get a 
decomposition 
\[ \phi : W_1 = V_{1,1}^* \oplus V_{1,2}^*
\xar{\left[ \begin{smallmatrix} \psi_{1,1}^* & 0 \\ 0 & 0 
\end{smallmatrix} \right] \cdot} 
V_{2,1}^* \oplus V_{2,2}^* = W_2  \]
in $\cat{TopMod} \k$. Now $W_1 \to V_{1,1}^*$ is a strict 
epimorphism, and $V_{1,1}^* \to W_2$ is a strict monomorphism. 

\medskip \noindent
(2) By the proof of part (1) we have isomorphisms
$\opn{Ker}(\phi) \cong V_{1,2}^*$, 
$\opn{Im}(\phi) \cong V_{1,1}^* \cong V_{2,1}^*$ and
$\opn{Coker}(\phi) \cong V_{2,2}^*$
in $\cat{TopMod} \k$. 

\medskip \noindent
(3) This is because
$\opn{Im}(\phi) = \opn{Ker}(W_2 \to V_{2,2}^*)$,
and $V_{2,2}^*$ is separated.
\end{proof}

\begin{cor} \label{cor4}
Let 
\[ S  = \left( 0 \to W_0 \xar{\phi_0} W_1 \xar{\phi_1} W_2 \to 0 
\right) \]
be a sequence of morphisms in $\cat{TopMod}_{\mrm{pf}} \k$. The 
following are equivalent:
\begin{enumerate}
\rmitem{i} The sequence $S$ is exact \tup{(}neglecting
topologies\tup{)}.
\rmitem{ii} The sequence $S$ is strict-exact in $\cat{TopMod} \k$.
\rmitem{iii} The sequence $S$ is split-exact in $\cat{TopMod} \k$.
\rmitem{iv} The dual sequence 
\[ S^*  := 
\left( 0 \to W_2^* \xar{\phi_1^*} W_1^* \xar{\phi_0^*} W_0^* \to 0 
\right) \]
 is exact in $\cat{Mod} \k$.
\end{enumerate}
\end{cor}

Condition (iii) says that there exists a continuous $\k$-linear 
homomorphism $\sigma : W_2 \to W_1$ such that
$\phi_1 \circ \sigma = \bsym{1}_{W_2}$; i.e.\ 
$W_1 \cong W_0 \oplus W_2$ in $\cat{TopMod} \k$. 

\begin{proof}
Any exact sequence in $\cat{Mod} \k$ splits; and by duality 
we deduce (iv) $\Rightarrow$ (iii). The implications (iii)
$\Rightarrow$ (ii) $\Rightarrow$ (i) are trivial. 
It remains to verify (i) $\Rightarrow$ (iv). 

Suppose $\psi : W \to W_1$ is some morphism in $\cat{TopMod} \k$ 
such that $\phi_1 \circ \psi = 0$. Because $\phi_0$ is a strict 
monomorphism it follows that $\psi$ factors through $\phi_0$. So 
$\phi_0$ is the kernel of $\phi_1$ in $\cat{TopMod} \k$ (in the 
categorical sense, see \cite[Sections II.6 and II.9]{HS}).

Next suppose $\psi : W_1 \to W$ is some morphism in $\cat{TopMod} \k$ 
such that $\psi \circ \phi_0 = 0$. Because $\phi_1$ is a strict 
epimorphism it follows that $\psi$ factors through $\phi_1$. So 
$\phi_1$ is the cokernel of $\phi_0$ in $\cat{TopMod} \k$ 

Since the categories $\cat{TopMod} \k$ and $\cat{Mod} \k$ are 
dual it follows that 
$\phi_1^* : W_2^* \to W_1^*$ is the kernel of 
$\phi_0^* : W_1^* \to W_0^*$, and vice versa. Therefore the 
sequence $S^*$ is exact.
\end{proof}

The classical notion of {\em linearly compact} topological 
$\k$-module is due to Lefchetz. 
See \cite[Sections 10.10 - 10.11]{Ko}.

\begin{cor}
Let $W$ be a topological $\k$-module. $W$ is profinite if and only 
if it is linearly compact.
\end{cor}

\begin{proof}
Say $W \in \cat{TopMod}_{\mrm{pf}} \k$. Let
$V := \opnt{D}^{\mrm{c}} W \in \cat{Mod} \k$.
By Theorem \ref{thm5}(2) we know that 
$W \cong \opnt{D} V$ in $\cat{TopMod} \k$.
According to Proposition \ref{prop1}(2), the topology on 
$\opnt{D} V$ is the weak$^*$ topology. Hence by 
\cite[Section 10.10 item (1)]{Ko} it follows that $\opnt{D} V$ 
is linearly compact.

Conversely, suppose $W$ is linearly compact. Let
$V$ be as above. By \cite[Section 10.10 item (3)]{Ko}  we get 
$W \cong \opnt{D} V$ in $\cat{TopMod} \k$,
and by Proposition \ref{prop1}(1) we see that 
$\opnt{D} V$ is profinite. 
\end{proof}

Given a set $X$ let
\[ \k^{(X)} := \{ \phi : X \to \k \mid \phi 
\text{ has finite support} \} . \]
This is a free $\k$-module with basis $\{ \delta_x \}_{x \in X}$,
where  $\delta_x : X \to \k$ is the ``delta function'' defined by 
\[ \delta_x(y) := \begin{cases}
1 & \text{if}\ y = x \\
0 & \text{otherwise} .
\end{cases} \]
If $V$ is a $\k$-module with basis $X$ then we get an isomorphism 
$\k^{(X)} \iso V$ by sending $\delta_x \mapsto x$. 

As usual $\k^X$ denotes the set of all functions 
$\phi : X \to \k$. We give it the product topology using the 
isomorphism $\k^X \cong \prod_{x \in X} \k$, where each copy of 
$\k$ has the discrete topology.

There are evaluation homomorphisms
$\alpha_X : \k^{(X)} \to \opnt{D}^{\mrm{c}}\, \k^{X}$ and
$\beta_X : \k^{X} \to \opnt{D}\, \k^{(X)}$
with formulas
\[ \beta_X(\phi)(\psi) = \alpha_X(\psi)(\phi)
:= \sum_{x \in X} \phi(x) \psi(x) \in \k  
\]
for $\phi \in \k^{X}$ and $\psi \in \k^{(X)}$.

Suppose $X$ and $Y$ are sets and $f : X \to Y$ is a function. 
Define a continuous $\k$-linear homomorphism 
$F^*(f) : \k^{Y} \to \k^{X}$
by 
\[ F^*(f)(\psi)(x) := \psi(f(x))  \]
for $\psi \in \k^Y$ and $x \in X$. Thus $F^*(f)(\psi)$ is the 
pullback of $\psi$. This is a functor
$F^* : \cat{Sets}^{\mrm{op}} \to \cat{TopMod} \k$.
Also define a $\k$-linear homomorphism 
$F(f) : \k^{(X)} \to \k^{(Y)}$
by 
\[ F(f)(\phi)(y) := \sum_{x \in f^{-1}(y)} \phi(x) \]
for $y \in Y$ and $\phi \in \k^{(X)}$. The function 
$F(f)(\phi)$ is the trace of $\phi$, or ``integration on the fibers 
of $f$''. We get a functor $F : \cat{Sets} \to \cat{Mod} \k$.

\begin{lem} \label{lem9}
Let $f : X \to Y$ be any function of sets. 
\begin{enumerate}
\item The diagram 
\[ \begin{CD}
\k^Y @>{F^*(f)}>> \k^X \\
@V{\beta_Y}VV @V{\beta_X}VV  \\
\opnt{D}\, \k^{(Y)} @>{\opnt{D}(F(f))}>> \opnt{D}\, \k^{(X)} 
\end{CD} \]
in $\cat{TopMod} \k$ is commutative.
\item The diagram 
\[ \begin{CD}
\k^{(X)} @>{F(f)}>> \k^{(Y)} \\
@V{\alpha_X}VV @V{\alpha_Y}VV  \\
\opnt{D}^{\mrm{c}}\, \k^{X} 
@>{\opnt{D}^{\mrm{c}}(F^*(f))}>> \opnt{D}^{\mrm{c}}\, \k^{Y}
\end{CD} \]
in $\cat{Mod} \k$ is commutative.
\end{enumerate}
\end{lem}

\begin{proof}
Direct calculation.
\end{proof}

\begin{prop} \label{prop5}
Let $X$ be any set. 
\begin{enumerate}
\item $\alpha_X : \k^{(X)} \to \opnt{D}^{\mrm{c}}\, \k^{X}$
is an isomorphism in $\cat{Mod} \k$.
\item $\beta_X : \k^{X} \to \opnt{D}\, \k^{(X)}$
is an isomorphism in $\cat{TopMod} \k$.
\end{enumerate}
\end{prop}

\begin{proof}
Denote by $\cat{Fin} X$ the set of finite subsets of $X$, 
which we rename $I$, and write $\cat{Fin} X = \{ X_i \}_{i \in I}$. 
$I$ is quasi-ordered by inclusion, and it is a directed set. 
The functor $F : \cat{Sets} \to \cat{Mod} \k$, 
$F Y = \k^{(Y)}$, makes
$\{ \k^{(X_i)} \}_{i \in I}$ into a direct system, and
$\k^{(X)} \cong \underset{i \in I}{\varinjlim}\, \k^{(X_i)}$
in $\cat{Mod} \k$. Moreover the set of finitely generated submodules
$\{ \k^{(X_i)} \}_{i \in I}$
is cofinal in $\cat{Fin} \k^{(X)}$. 

The opposite quasi-ordered set $I^{\mrm{op}}$ is also directed. 
The functor $F^* : \cat{Sets}^{\mrm{op}} \to \cat{TopMod} \k$, 
$F^* Y = \k^{Y}$, makes
$\{ \k^{X_i} \}_{i \in I}$ into an inverse system, and
$\k^{X} \cong \underset{i \in I}{\varprojlim}\, \k^{X_i}$
in $\cat{TopMod} \k$. The inverse system of open cofinite submodules
$\{ \opn{Ker}(\k^{X} \to \k^{X_i}) \}_{i \in I^{\mrm{op}}}$
is cofinal in $(\cat{Cofin} \k^X)^{\mrm{op}}$. 

By Proposition \ref{prop4}, $\k^X$ is a profinite topological 
$\k$-module; hence by Lemma \ref{lem4} we have isomorphisms
\[ \opnt{D}^{\mrm{c}}\, \k^X = 
\opn{Hom}_{\k}^{\mrm{cont}}(\k^X, \k) \cong
\underset{W' \in \cat{Cofin} \k^X}{\varinjlim}\, 
\opnt{D}^{\mrm{c}}\, (\k^X / W') 
\cong \underset{i \in I}{\varinjlim}\, \opnt{D}^{\mrm{c}}\, \k^{X_i}
. \]
Now for the finite set $X_i$ one has
$\k^{(X_i)} = \k^{X_i}$, and the adjunction map
$\alpha_{X_i} : \k^{X_i} \to \opnt{D}^{\mrm{c}}\, \k^{X_i}$
is bijective. So by Lemma \ref{lem9}(2) we can switch the limits 
from the functor $F^*$ to the functor $F$:
\[ \underset{i \in I}{\varinjlim}\, \opnt{D}^{\mrm{c}}\, \k^{X_i}
\cong
\underset{i \in I}{\varinjlim}\, \k^{(X_i)} \cong  \k^{(X)} . \]
We get an isomorphism $\opnt{D}^{\mrm{c}}\, \k^X  \cong \k^{(X)}$ 
that is compatible with the evaluation pairing, so $\alpha_X$ is an 
isomorphism. 

Finally we look at $\beta_X$. We know that 
$\{ \k^{(X_i)} \}_{i \in I}$
is cofinal in $\cat{Fin} \k^{(X)}$, so 
\[ \opnt{D}\, \k^{(X)} = 
\underset{V' \in \cat{Fin} \k^{(X)}}{\varprojlim}\, \opnt{D}\, V'
\cong 
\underset{i \in I}{\varprojlim}\, \opnt{D}\, \k^{(X_i)} \]
in $\cat{TopMod} \k$. But using the isomorphisms
$\beta_{X_i} : \k^{X_i} \iso \opnt{D}\, \k^{(X_i)}$
and Lemma \ref{lem9}(1) we can switch 
limits to obtain an isomorphism
\[ \underset{i \in I}{\varprojlim}\, \opnt{D}\, \k^{(X_i)} \cong
\underset{i \in I}{\varprojlim}\, \k^{X_i} \cong \k^X \]
in $\cat{TopMod} \k$. These isomorphisms are compatible with 
evaluation. Therefore $\beta_X$ is an isomorphism.
\end{proof}

\begin{cor} \label{cor9}
Let $W \in \cat{TopMod}_{\mrm{pf}} \k$. Then
$W \cong \k^X$ in $\cat{TopMod} \k$ for some set $X$.
\end{cor}

\begin{proof}
Choose a basis $X$ for $\opnt{D}^{\mrm{c}} W$. 
Then $\opnt{D}^{\mrm{c}} W \cong \k^{(X)}$ in 
$\cat{Mod} \k$, and according to Proposition \ref{prop5} and 
Theorem \ref{thm5} we get
\[ W \cong \opnt{D} \opnt{D}^{\mrm{c}} W \cong
\opnt{D}\, \k^{(X)} \cong \k^{X} . \]
\end{proof}

\begin{rem}
Not all cofinite submodules of a profinite topological $\k$-module 
are open. Take $\k := \mbb{Q}$ and $W := \k^{\mbb{N}}$. Since the 
open cofinite submodules of $W$ are parameterized by the 
submodules of $\k^{ \{ 0, \ldots, n \} }$ for $n \in \mbb{N}$ it 
follows that $\lvert \cat{Cofin} W \rvert = \aleph_0$. 
But on the other 
hand $\opn{rank}_{\k} W = \lvert W \rvert = 2^{\aleph_0}$, 
so there are at least this 
many cofinite submodules $W' \subset W$. 
\end{rem}

\begin{rem}
It is curious to note that some $\k$-modules $W$ do not admit any 
topology with which they are profinite. Indeed, when 
$\k = \mbb{Q}$ the module $W := \k^{(\mbb{N})}$ is such an example. 
Its rank is $\aleph_0$. Now for a $\k$-module $V$ one has either
$\opn{rank}_{\k} V < \aleph_0$, and then
$\opn{rank}_{\k} \opnt{D} V = \opn{rank}_{\k} V < \aleph_0$; or 
$\opn{rank}_{\k} V \geq \aleph_0$, in which case
$\opn{rank}_{\k} \opnt{D} V \geq 2^{\aleph_0}$. 
So $W \cong \opnt{D} V$ is impossible for any $V$.
\end{rem}

\begin{rem}
The duality here is really the Lefchetz duality for linearly 
compact topological vector spaces. See \cite[Section 10.10]{Ko}.
The way it is presented here lends itself easily to 
generalizations. Indeed, all definitions and results in this 
section, up to and including Corollary \ref{cor8}, are valid for 
any artinian commutative local ring $\k$. 
The sole modification needed 
is that the dualities should be
$\opnt{D} := \opn{Hom}_{\k}(-, J)$
and
$\opnt{D}^{\mrm{c}} := \opn{Hom}^{\mrm{cont}}_{\k}(-, J)$,
where $J$ is some fixed injective hull of the residue field of 
$\k$. Likewise all results on function modules (except 
Corollary \ref{cor9}) hold, once we 
replace $\k^{(X)}$ with $J^{(X)}$ everywhere.

In even greater generality we arrive at 
the Gabriel-Matlis theory of pseudo-compact modules. Another 
variant is the duality theory of Beilinson completion algebras
in \cite{Ye2}. 
\end{rem}

\section{Topological $A$-Modules}

Now let $A$ be an associative unital algebra over the field $\k$.
We do not assume $A$ is commutative nor (left or right) 
noetherian.
All $A$-modules are by default left modules. We have the category 
$\cat{TopMod} A$ consisting of topological $A$-modules (any sort 
of topology) and continuous $A$-linear homomorphisms. Like 
$\cat{TopMod} \k$, the category $\cat{TopMod} A$ is additive, and 
it has exact sequences (in which the homomorphisms are 
required to be strict). There are forgetful functors
$\cat{TopMod} A \to \cat{TopMod} \k$ and
$\cat{TopMod} A \to \cat{Mod} A$. 

\begin{dfn} \label{dfn2}
A topological $A$-module $M$ is called {\em profinite over $\k$} 
if it is profinite when considered as topological $\k$-module. 
We denote by $\cat{TopMod}_{\mrm{pf} / \k} A$ the full subcategory 
of $\cat{TopMod} A$ consisting of topological $A$-modules 
profinite over $\k$.
\end{dfn}

We denote by $A^{\mrm{op}}$ the opposite algebra, i.e.\ the same 
$\k$-module but with reversed multiplication. A right $A$-module 
is then the same as a left $A^{\mrm{op}}$-module. 
When $A$ is commutative of course $A^{\mrm{op}} = A$.

\begin{thm} \label{thm6}
Let $A$ be a $\k$-algebra.
Given $M \in \cat{Mod} A$ the dual $M^* = \opnt{D} M$ is a 
topological $A^{\mrm{op}}$-module profinite over $\k$. The functor 
\[ \opnt{D} : \cat{Mod} A \to 
\cat{TopMod}_{\mrm{pf} / \k} A^{\mrm{op}} \]
is an equivalence, with adjoint $\opnt{D}^{\mrm{c}}$.
\end{thm}

\begin{proof}
The $A^{\mrm{op}}$-module structure on $\opnt{D} M$ 
is defined as follows. 
Given an element $a \in A$ the function $\phi_a : M \to M$, 
$\phi_a(m) := a m$, is $\k$-linear, and so 
we get a continuous $\k$-linear homomorphism
$\opnt{D}(\phi_a) : \opnt{D} M \to \opnt{D} M$. Because $\opnt{D}$ 
is a contravariant functor we have
\[ \opnt{D}(\phi_{a_1 a_2}) = 
\opnt{D}(\phi_{a_1} \circ \phi_{a_2}) = 
\opnt{D}(\phi_{a_2}) \circ \opnt{D}(\phi_{a_1})  \]
and
\[ \opnt{D}(\phi_{1}) = \opnt{D}(\bsym{1}_M) = 
\bsym{1}_{\opnt{D} M} \]
for the element $1 \in A$. Likewise we get an $A$-module structure 
on $\opnt{D}^{\mrm{c}} N$ for any 
$N \in \cat{TopMod}_{\mrm{pf} / \k} A^{\mrm{op}}$. 
One sees easily that the adjunction homomorphism
$\beta_N : N \to \opnt{D} \opnt{D}^{\mrm{c}} N$ 
is $A^{\mrm{op}}$-linear, and
$\alpha_M : M \to \opnt{D}^{\mrm{c}} \opnt{D} M$ 
is $A$-linear. According to Theorem \ref{thm5} the homomorphisms
$\beta_N$ and $\alpha_M$ are isomorphisms in the respective 
categories.
\end{proof}

\begin{cor} \label{cor5}
Let 
\[ S  = \left( \cdots \to N_0 \xar{\phi_0} N_1 \xar{\phi_1} 
N_2 \to \cdots \right) \]
be a sequence of morphisms in $\cat{TopMod}_{\mrm{pf} / \k} A$. The 
following are equivalent:
\begin{enumerate}
\rmitem{i} The sequence $S$ is exact \tup{(}neglecting
topologies\tup{)}.
\rmitem{ii} The sequence $S$ is strict-exact.
\rmitem{iii} The dual sequence 
\[ S^*  := 
\left( \cdots \to N_2^* \xar{\phi_1^*} N_1^* \xar{\phi_0^*} 
N_0^* \to \cdots \right) \]
is exact in $\cat{Mod} A^{\mrm{op}}$.
\end{enumerate}
\end{cor}

\begin{proof}
Combine Theorem \ref{thm6} and Corollary \ref{cor4}.
\end{proof}

\begin{cor} \label{cor6}
Let $N$ be a topological $A$-module. Then the following are 
equivalent:
\begin{enumerate}
\rmitem{i} $N$ is isomorphic in $\cat{TopMod} A$ to a closed 
submodule of $(A^*)^r = (\opnt{D} A)^r$ for some natural number 
$r$. 
\rmitem{ii} $N$ is a profinite topological $\k$-module and 
$N^* = \opnt{D}^{\mrm{c}} N$ is a finitely generated 
$A^{\mrm{op}}$-module. 
\end{enumerate}
\end{cor}

\begin{proof}
(i) $\Rightarrow$ (ii): Let $\psi : N \to (A^*)^r$ be a
continuous $A$-linear homomorphism such that $\opn{Im}(\psi)$ is a 
closed submodule of $(A^*)^r$ and $\psi : N \to \opn{Im}(\psi)$ is 
an isomorphism. By Theorem \ref{thm0} $\opn{Im}(\psi)$ is a 
profinite topological $\k$-module, and hence so is $N$. The 
homomorphism $\psi : N \to (A^*)^r$ is a 
monomorphism in $\cat{TopMod}_{\mrm{pf} / \k} A$, so according to 
Corollary \ref{cor5} the dual $\psi^* : A^r \to N^*$
is surjective. Thus $N^*$ is a finitely generated 
$A^{\mrm{op}}$-module. 

\medskip \noindent
(ii) $\Rightarrow$ (i): Let $M := N^*$. There is an epimorphism
$\phi : A^r \to M$ for some $r$. Dualizing, and using Corollary 
\ref{cor5}, we get a monomorphism $\phi^* : M^* \to (A^*)^r$. 
Because $N$ is profinite over $\k$ there is an isomorphism
$N \cong M^*$ in $\cat{TopMod}_{\mrm{pf} / \k} A$. 
\end{proof}

Let $A$ be any $\k$-algebra. Then $A^* = \mrm{D} A$ is a 
topological $A$-bimodule, i.e.\ a topological module over 
$A \otimes_{\k} A^{\mrm{op}}$. For a natural number $r$ we 
consider elements of $(A^*)^r$ as rows of size $r$. Thus given 
$r_0, r_1 \in \mbb{N}$ and a matrix
$G \in \mrm{M}_{r_0 \times r_1}(A)$,
right multiplication by $G$ is a continuous $A$-linear 
homomorphism
$(A^*)^{r_0} \xar{\cdot G} (A^*)^{r_1}$.

\begin{cor} \label{cor2}
Suppose $A$ is a noetherian $\k$-algebra,
$M \in \cat{TopMod}_{\mrm{pf} / \k} A$, $r_0 \in \mbb{N}$,
and $\psi : M \to (A^*)^{r_0}$ is an
injective continuous $A$-linear homomorphism.
Then there exists some $r_1 \in \mbb{N}$, a matrix 
$G \in \mrm{M}_{r_0 \times r_1}(A)$, and an exact sequence
\begin{equation} \label{eqn4}
0 \to M \xar{\psi} (A^*)^{r_0} \xar{\cdot G} 
(A^*)^{r_1}
\end{equation}
in $\cat{TopMod} A$.
\end{cor}

\begin{proof}
By Corollary \ref{cor5} we get an epimorphism
$\psi^* : A^{r_0} \to M^*$ in $\cat{Mod} A^{\mrm{op}}$. 
Because $A^{\mrm{op}}$ is noetherian the kernel of $\psi^*$ 
is a finitely generated $A^{\mrm{op}}$-module, and therefore there 
is an exact sequence
\begin{equation} \label{eqn2}
A^{r_1} \to A^{r_0} \xar{\psi^*} M^* \to 0
\end{equation}
in $\cat{Mod} A^{\mrm{op}}$. If we think of $A^{r_i}$ as column 
vectors then the homomorphism $A^{r_1} \to A^{r_0}$ is given by left 
multiplication with some matrix 
$G \in \mrm{M}_{r_0 \times r_1}(A)$.
To finish we apply the duality functor $\opnt{D}$ to the sequence
(\ref{eqn2}).
\end{proof}

\begin{dfn}
Let $M \in \cat{TopMod}_{\mrm{pf} / \k} A$. An exact sequence 
\[ 0 \to M \to (A^*)^{r_0} \xar{\cdot G} (A^*)^{r_1} \]
in $\cat{TopMod} A$ (cf.\ Corollary \ref{cor5}), 
for some $r_0, r_1 \in \mbb{N}$ and some matrix 
$G \in$ \linebreak $\mrm{M}_{r_0 \times r_1}(A)$, is called a 
{\em kernel representation} of $M$. 
\end{dfn}

\begin{cor} \label{cor3}
Let $A$ be a commutative $\k$-algebra, let
$M, N \in \cat{TopMod}_{\mrm{pf} / \k} A$, and suppose we are 
given kernel representations
\[ 0 \to M \to (A^*)^{r_0} \xar{\cdot G} (A^*)^{r_1} 
\]
and
\[ 0 \to N \to (A^*)^{s_0} \xar{\cdot F} (A^*)^{s_1} ,
\]
where $G \in \mrm{M}_{r_0 \times r_1}(A)$
and $F \in \mrm{M}_{s_0 \times s_1}(A)$. Let
$\phi : M \to N$ be a continuous $A$-linear homomorphism. 
Then there exist matrices
$H_0$ and $H_1$ of appropriate sizes with 
entries in $A$, such that the diagram 
\[ \begin{CD}
0 @>>> M @>>> (A^*)^{r_0} @>{\cdot G}>> 
(A^*)^{r_1} \\
& &  @V{\phi}VV @VV{\cdot H_0}V @VV{\cdot H_1}V  \\
0 @>>> N @>>> (A^*)^{s_0} @>{\cdot F}>> 
(A^*)^{s_1} 
\end{CD} \]
is commutative.
\end{cor}

\begin{proof}
By duality we obtain exact sequences
\[ A^{r_1} \xar{G \cdot} A^{r_0} \to M^* \to 0 \]
and
\[ A^{s_1} \xar{F \cdot} A^{s_0} \to N^* \to 0  \]
in $\cat{Mod} A$. Because the modules $A^{r_i}$ and $A^{s_i}$ 
are free the homomorphism 
$\phi^* : N^* \to M^*$ extends to a commutative diagram
\[ \begin{CD}
A^{s_1} @>{F \cdot}>> A^{s_0} @>>> N^* @>>> 0 \\
@VV{}V @VV{}V @V{\phi^*}VV  \\
A^{r_1} @>{G \cdot}>> A^{r_0} @>>> M^* @>>> 0 
\end{CD} \]
with exact rows (see \cite[Theorem 6.9]{Ro}). The homomorphism
$A^{s_i} \to A^{r_i}$ is left multiplication by some matrix
$H_i \in \mrm{M}_{r_i \times s_i}(A)$. 
Now apply the functor $\opnt{D}$ to this diagram. 
\end{proof}

\section{Behaviors}

Let $n$ be a positive integer. For any $1 \leq j \leq n$
the shift $s_j : \mbb{N}^n \to \mbb{N}^n$  
is defined by
\[ s_j(i_1, \ldots, i_n) := (i_1, \ldots, i_j + 1,  \ldots, i_n) 
. \]
Recall that $\k$ is a field and
$\k^{\mbb{N}^n} = \{ \phi : \mbb{N}^n \to \k \}$. 

\begin{dfn} \label{dfn3}
The {\em $n$-dimensional shift module} is the 
topological $\k$-module  
$\k^{\mbb{N}^n} = \prod_{\bsym{i} \in \mbb{N}^n} \k$,
with the product topology. The shift operators 
$\sigma_1, \ldots, \sigma_n : \k^{\mbb{N}^n} \to \k^{\mbb{N}^n}$ 
are defined as follows:
$\sigma_j(\phi) := \phi \circ s_j$
for $\phi \in \k^{\mbb{N}^n}$.
\end{dfn}

\begin{dfn} \label{dfn4}
Let $n$ be a positive integer. 
A {\em $n$-dimensional behavior} is a closed $\k$-submodule 
$M \subset (\k^{\mbb{N}^n})^{r}$
for some nonnegative integer $r$ that is invariant under the shift 
operators $\sigma_1, \ldots, \sigma_n$.
\end{dfn}

Let $\k[\bsym{z}] = \k[z_1, \ldots, z_n]$ be the polynomial 
algebra in $n$ indeterminates. Since the shift operators commute, 
the shift module $\k^{\mbb{N}^n}$ is a topological 
$\k[\bsym{z}]$-module, where $z_j$ acts by the operator 
$\sigma_j$. So an $n$-dimensional behavior is precisely a closed 
$\k[\bsym{z}]$-submodule of $(\k^{\mbb{N}^n})^{r}$.

\begin{prop} \label{prop6}
One has 
$\k^{\mbb{N}^n} \cong \opnt{D}\, \k[\bsym{z}] = \k[\bsym{z}]^*$ 
in $\cat{TopMod} \k[\bsym{z}]$. Thus in particular
$\k^{\mbb{N}^n} \in \cat{TopMod}_{\mrm{pf} / \k} \k[\bsym{z}]$. 
\end{prop}

\begin{proof}
The set of monomials 
$\{ \bsym{z}^{\bsym{i}} \}_{\bsym{i} \in \mbb{N}^n} \cong \mbb{N}^n$
is a basis of the $\k$-module $\k[\bsym{z}]$, and checking
the shift action we see that
$\k[\bsym{z}] \cong \k^{(\mbb{N}^n)}$ as $\k[\bsym{z}]$-modules.
Now use Proposition \ref{prop5}(2).
\end{proof}

Let $\k[[x]]$ be the ring of formal power series and 
$\k((x))$ the field of Laurent series, with their usual topologies 
(cf.\ \cite[Section 1.3]{Ye1}). Let us write $x := z^{-1}$. 
Since $\k[z] \subset \k((x)) = \k((z^{-1}))$ 
is a subring we get an exact sequence of $\k[z]$-modules
\begin{equation} \label{eqn3}
0 \to \k[z] \to \k((z^{-1})) \to 
z^{-1}\, \k[[z^{-1}]] \to 0 .
\end{equation}
Put the discrete topology on $\k[z]$ and the product topology on 
$z^{-1}\, \k[[z^{-1}]] \cong \prod_{i \leq -1} \k$. 
Then the sequence (\ref{eqn3}) is split exact in 
$\cat{TopMod} \k$, see \cite[Proposition 1.3.5]{Ye1}. 
(We shall not need this fact.) The next lemma is clear. 

\begin{lem} \label{lem3}
The map
$z^{-1}\, \k[[z^{-1}]] \to \k^{\mbb{N}}$
sending
\[ \sum_{i \in \mbb{N}} \lambda_i z^{-1 -i} \mapsto 
(\lambda_0, \lambda_1, \ldots) \]
is an isomorphism in $\cat{TopMod} \k[z]$. 
\end{lem}

Traditionally the shift module is defined as $z^{-1}\, \k[[z^{-1}]]$,
cf.\ \cite{Fu2}. In view of Lemma \ref{lem3} a behavior in the 
sense of \cite{Fu2} is precisely a $1$-dimensional behavior in 
the sense of Definition \ref{dfn4}.

\begin{rem} 
Not all $\k[z]$-submodules of $\k^{\mbb{N}}$ are closed. 
Take 
\[ w := (1, 0, 1, 0, 0, 1, 0, 0, 0, 1, \ldots) . \]
Then the submodule 
$N := \k[z] \, w \subset \k^{\mbb{N}}$ is dense in 
$\k^{\mbb{N}}$ but not equal to it.
\end{rem}

Here is a classification of behaviors.

\begin{thm} \label{thm2}
Let $M$ be an $n$-dimensional behavior and let
$\psi : M \to (\k^{\mbb{N}^n})^{r_0}$ be an injective continuous 
$\k[\bsym{z}]$-linear homomorphism. 
\begin{enumerate}
\item There exist natural numbers
$r_1, \ldots, r_n$ and matrices
$G_i(\bsym{z}) \in \mrm{M}_{r_{i - 1} \times r_{i}}(\k[\bsym{z}])$ 
such that the sequence of homomorphisms
\[ 0 \to M \xar{\psi} (\k^{\mbb{N}^n})^{r_0} 
\xar{\cdot G_1(\bsym{z})} 
(\k^{\mbb{N}^n})^{r_1} \xar{\cdot G_2(\bsym{z})} 
\cdots \xar{\cdot G_n(\bsym{z})} (\k^{\mbb{N}^n})^{r_n} \to 0 \]
is exact, and in fact splits in $\cat{TopMod} \k$. 
\item If $n = 1$ there is an isomorphism
\[ M \cong (\k^{\mbb{N}})^r \oplus N \]
of topological $\k[z]$-modules, where $r := r_0 - r_1$ and $N$ is 
finitely generated as $\k$-module.
\end{enumerate}
\end{thm}

\begin{proof}
(1) From Corollary \ref{cor2} we get an epimorphism 
$\psi^* : \k[\bsym{z}]^{r_0} \to M^*$
in $\cat{Mod} \k[\bsym{z}]$.
The Hilbert Syzygy Theorem (cf.\ \cite[Corollary 9.3.5]{Ro})
says we can extend it to an exact sequence
\begin{equation} \label{eqn1}
0 \to \k[\bsym{z}]^{r_n} \xar{G_n(\bsym{z}) \cdot}
\cdots \xar{G_2(\bsym{z}) \cdot}
\k[\bsym{z}]^{r_1} \xar{G_1(\bsym{z}) \cdot} \k[\bsym{z}]^{r_0} 
\xar{\psi^*} M^* \to 0
\end{equation}
for some $r_i$ and some matrices $G_i(\bsym{z})$. This sequence 
splits in $\cat{Mod} \k$. 
Now apply the functor $\opnt{D}$. 

\medskip \noindent 
(2) In the case $n = 1$, by the theory of finitely 
generated modules over a 
PID we know that 
$M^* \cong \k[z]^r \oplus L$ for some $r$ and some torsion 
$\k[z]$-module $L$. Moreover tensoring the sequence (\ref{eqn1}) 
with the field $\k(z)$
we see that $r = r_0 - r_1$. Applying the functor $\opnt{D}$ to the 
isomorphism $M^* \cong \k[z]^r \oplus L$ we get
$M \cong (\k^{\mbb{N}})^r \oplus L^*$.
\end{proof}

\begin{thm} \label{thm4}
Let $M$ and $N$ be two $n$-dimensional
behaviors, and let $\phi : M \to N$ be a 
continuous $\k[\bsym{z}]$-linear homomorphism. 
Suppose we are given kernel representations
\[ 0 \to M \to (\k^{\mbb{N}^n})^{r_0} \xar{\cdot G_1(\bsym{z})} 
(\k^{\mbb{N}^n})^{r_1} \]
and
\[ 0 \to N \to (\k^{\mbb{N}^n})^{s_0} \xar{\cdot F_1(\bsym{z})} 
(\k^{\mbb{N}^n})^{s_1} \]
of these behaviors. Then there exist matrices
$H_0(\bsym{z})$ and $H_1(\bsym{z})$ of appropriate sizes with 
entries in $\k[\bsym{z}]$ such that the diagram 
\[ \begin{CD}
0 @>>> M @>>> (\k^{\mbb{N}^n})^{r_0} @>{\cdot G_1(\bsym{z})}>> 
(\k^{\mbb{N}^n})^{r_1} \\
& &  @V{\phi}VV @VV{\cdot H_0(\bsym{z})}V @VV{\cdot H_1(\bsym{z})}V 
\\ 
0 @>>> N @>>> (\k^{\mbb{N}^n})^{s_0} @>{\cdot F_1(\bsym{z})}>> 
(\k^{\mbb{N}^n})^{s_1} 
\end{CD} \]
is commutative.
\end{thm}

When $n = 1$ this is \cite[Theorem 3.4]{Fu2}. 

\begin{proof}
This is a special case of Corollary 
\ref{cor3} for the ring $A := \k[\bsym{z}]$.
\end{proof}

\section{Noncontinuous Homomorphisms}

In this section $A$ is a finitely generated commutative 
algebra over the field $\k$. We analyze the algebraic structure 
of the $A$-module $A^*$, namely we forget the 
topology. This part requires more difficult ring theory.

As usual we shall denote by $\opn{Spec} A$ the set of prime 
ideals of the ring $A$, and by $\opn{Max} A$ the subset of maximal 
ideals. For any prime ideal $\mfrak{p}$
let $J(\mfrak{p})$ be an injective hull of 
$A / \mfrak{p}$ considered as $A$-module. 
Cf.\ \cite[Section 3]{Ro} or \cite[Section V.4]{St}. 

\begin{prop} \label{prop7}
There is a \tup{(}noncanonical\tup{)}
decomposition of $A$-modules
\[ A^* \cong 
\Bigl( \bigoplus_{\mfrak{m} \in \opn{Max} A}
J(\mfrak{m}) \Bigr) \oplus N . \]
Here $N$ is some injective $A$-module that contains no nonzero 
finite length submodules.
\end{prop}

\begin{proof}
Since for any $A$-module $M$ there is a functorial isomorphism 
\[ \opn{Hom}_{A}(M, A^*) \cong M^* \]
it follows that $A^*$ is an injective $A$-module. 
Because $A$ is a noetherian ring we know that 
\[ A^* \cong \bigoplus_{\mfrak{p} \in \opn{Spec} A}
J(\mfrak{p})^{(\mu_{\mfrak{p}})} , \]
where each $\mu_{\mfrak{p}}$ is a cardinal number, and
$J(\mfrak{p})^{(\mu_{\mfrak{p}})}$
is a direct sum of $\mu_{\mfrak{p}}$ copies of $J(\mfrak{p})$.
See \cite[Propositions V.4.5 and V.4.6]{St}.

Let $\mfrak{m}$ be some maximal ideal and let
$\bsym{k}(\mfrak{m}) := A / \mfrak{m}$, the residue field.
Then 
\[ \mu_{\mfrak{m}} = \opn{rank}_{\bsym{k}(\mfrak{m})}
\opn{Hom}_A(\bsym{k}(\mfrak{m}), A^*) . \]
On the other hand
\[ \opn{Hom}_A(\bsym{k}(\mfrak{m}), A^*) \cong 
\bsym{k}(\mfrak{m})^* \cong  \bsym{k}(\mfrak{m}) \]
as $\bsym{k}(\mfrak{m})$-modules. Hence 
$\mu_{\mfrak{m}} = 1$. 

Finally 
$N := \bigoplus_{\mfrak{p}} J(\mfrak{p})^{(\mu_{\mfrak{p}})}$,
the sum going over all prime ideals that are not maximal.
\end{proof}

\begin{cor} \label{cor7}
Let $n$ be any positive integer and 
$\k[\bsym{z}] := \k[z_1, \ldots, z_n]$. 
There exist $\k[\bsym{z}]$-linear homomorphisms
$\phi : \k^{\mbb{N}^n} \to \k^{\mbb{N}^n}$ that are not 
continuous. 
\end{cor}

\begin{proof}
 From Theorem \ref{thm6} and Proposition \ref{prop6} we know that
\[ \opn{Hom}^{\mrm{cont}}_{\k[\bsym{z}]}(\k^{\mbb{N}^n}, 
\k^{\mbb{N}^n}) \cong
\opn{Hom}_{\k[\bsym{z}]}(\k[\bsym{z}], \k[\bsym{z}]) \cong
\k[\bsym{z}] . \]
Thus a continuous homomorphism 
$\phi : \k^{\mbb{N}^n} \to \k^{\mbb{N}^n}$ 
is multiplication by some 
polynomial $f(\bsym{z}) \in \k[\bsym{z}]$. 
It follows that $\phi$ must preserve the decomposition in 
Proposition \ref{prop7}. Moreover, because
\[ \opn{Hom}_{\k[\bsym{z}]}(J(\mfrak{m}), J(\mfrak{m})) 
\cong \widehat{\k[\bsym{z}]}_{\mfrak{m}} , \]
where the latter is the $\mfrak{m}$-adic completion of 
$\k[\bsym{z}]$, we see that $\phi = 0$ iff the restriction
$\phi|_{J(\mfrak{m})} = 0$, where $\mfrak{m}$ is any maximal 
ideal.

We now go about constructing our counterexample. Choose a particular 
maximal ideal $\mfrak{m}_0$; say 
$\mfrak{m}_0 := (z_1, \ldots, z_n)$. Using the 
decomposition of Proposition \ref{prop7} define 
\[ \phi|_{J(\mfrak{m}_0)} : J(\mfrak{m}_0) \to \k^{\mbb{N}^n} \]
to be the inclusion; define
\[ \phi|_{J(\mfrak{m})} : J(\mfrak{m}) \to \k^{\mbb{N}^n} \]
to be the zero homomorphism for any other maximal ideal $\mfrak{m}$;
and define
\[ \phi|_{N} : N \to \k^{\mbb{N}^n} \]
also to be the zero homomorphism. This homomorphism cannot be 
continuous.
\end{proof}

\begin{rem}
When $A := \k[z]$, i.e.\ the one dimensional case, the module $N$ 
appearing in Proposition \ref{prop7} is of the form
$N \cong \k(z)^{(\mu_0)}$. Regarding the cardinality $\mu_0$, we 
know it if $\k$ is a countable field. In this case we have
\[ \mu_0 = \opn{rank}_{\k(z)} N = \abs{N} = 
\abs{\k^{\mbb{N}}} = 2^{\aleph_0} = \aleph . \]
\end{rem}

\section{Noncommutative $n$-Dimensional Behaviors}

Fix a positive integer $n$. Let $S$ be the free monoid 
(i.e.\ a semigroup with $1$) on the generators $s_1, \ldots, s_n$. 
So the elements of $S$ are the words $s_{i_1} \cdots s_{i_k}$
with $k \in \mbb{N}$ and 
$i_1, \ldots, i_k \in \{ 1, \ldots, n \}$. 

Suppose $X$ is a set equipped with functions
$f_1, \ldots, f_n : X \to X$. This data determines a left action 
of the monoid $S$ on $X$, by letting
$s_i(x) := f_i(x)$ for the generators $s_1, \ldots, s_n$ of $S$.
A set $X$ with a left action of $S$ is called a {\em left 
$S$-set}. 

\begin{dfn}
A {\em finitely generated left $S$-set} is a left $S$-set $X$, 
with a finite subset $X_0 \subset X$, such that 
$X = \{ s(x) \mid s \in S \text{ and } x \in X_0 \}$. 
\end{dfn}

Recall that $\k$ is a field. Given a left $S$-set $X$ there is a 
right action of $S$ on $\k^X$ by continuous $\k$-linear 
homomorphisms, which we call {\em shift operators}. 

\begin{dfn}
A {\em noncommutative $n$-dimensional behavior} is a closed shift 
invariant $\k$-submodule of $\k^X$, for some finitely generated 
left $S$-set $X$.
\end{dfn}

Let $\k \bra{\bsym{z}} = \k \bra{z_1, \ldots, z_n}$
be the free associative algebra on the variables 
$z_1, \ldots,$ \linebreak $z_n$. 
We can identify the monoid $S$ with the 
multiplicative monoid of monomials in $\k \bra{\bsym{z}}$.
In this way any $n$-dimensional behavior $M$ becomes a right 
$\k \bra{\bsym{z}}$-module, i.e.\ 
$M \in \cat{TopMod}_{\mrm{pf} / \k} 
\k \bra{\bsym{z}}^{\mrm{op}}$.

If $r$ is an infinite cardinal number and $M$ is a $\k$-module,
then $M^{(r)}$ denotes direct sum of $r$ copies of $M$, 
whereas $M^r$ is the direct product of $r$ copies of $M$.

Recall that for $M \in \cat{TopMod} \k$ we write
$M^* := \mrm{D}^{\mrm{c}} M = \opn{Hom}_{\k}^{\mrm{cont}}
(M, \k)$; 
whereas for $M \in \cat{Mod} \k$ we write
$M^* := \mrm{D} M = \opn{Hom}_{\k}(M, \k)$.

\begin{thm}
Let $M$ be a noncommutative $n$-dimensional behavior. Then there 
exist $r_0 \in \mbb{N}$, $r_1 \in \mbb{N} \cup \{ \aleph_0 \}$,
$G(\bsym{z}) \in \mrm{M}_{r_0 \times r_1}(\k \bra{\bsym{z}})$ and 
a homomorphism 
$\phi : M \to (\k^S)^{r_0}$ in 
$\cat{TopMod} \k \bra{\bsym{z}}^{\mrm{op}}$
such that 
\begin{equation} \label{eqn5}
0 \to M \xar{\phi} (\k^{S})^{r_0} \xar{\cdot G(\bsym{z})} 
(\k^{S})^{r_1} \to 0
\end{equation}
is an exact sequence.
\end{thm}

\begin{proof}
By definition $M \subset \k^X$ for some finitely generated $S$-set 
$X$. Dualizing we obtain a surjection 
$\k^{(X)} \surj M^*$ in $\cat{Mod} \k \bra{\bsym{z}}$.
Since $X$ is finitely generated there is a surjection of sets 
$\underset{r_0}{\underbrace{S \sqcup \cdots \sqcup S}}
\surj X$
for some $r_0 \in \mbb{N}$. This gives rise to a surjection
$(\k^{(S)})^{r_0} \surj \k^{(X)}$. Composing we obtain a 
$\k \bra{\bsym{z}}$-linear surjection
$\phi^* : (\k^{(S)})^{r_0} \surj M^*$. 

Now $(\k^{(S)})^{r_0}$ is a free $\k \bra{\bsym{z}}$-module of 
finite rank, so according to 
\cite[Section 1.2 Theorem 2.1 and Section 2.4 Corollary 4.3]{Co}, 
the submodule $\opn{Ker}(\phi^*)$ is free. Let  
$r_1 \in \mbb{N} \cup \{ \aleph_0 \}$
be the rank of $\opn{Ker}(\phi^*)$. Hence there is an 
exact sequence
\[ 0 \to (\k^{(S)})^{(r_1)} \to (\k^{(S)})^{r_0} \to 
M^* \to 0 \]
in $\cat{Mod} \k \bra{\bsym{z}}$. Finally apply duality.
\end{proof}

\begin{exa}
Let $M$ be an $n$-dimensional (commutative) behavior. By definition 
$M$ is a closed shift invariant $\k$-submodule of $\k^X$, where
$X := 
\underset{r_0}{\underbrace{\mbb{N}^n \sqcup \cdots \sqcup 
\mbb{N}^n}}$. We see that $M$ is also a noncommutative 
$n$-dimensional behavior.
\end{exa}

Sometimes the exponent $r_1$ in a kernel representation of a 
noncommutative behavior must be infinite. Here is an example. 

\begin{exa}
Take $n = 2$, and let $M := \k^{\mbb{N}^2}$. The dual module 
$M^*$ is isomorphic, as $\k \bra{z_1, z_2}$-module, 
to the commutative polynomial ring $\k[z_1, z_2]$. Consider the 
surjection
$\phi^* : \k \bra{z_1, z_2} \surj \k[z_1, z_2]$. 
Then $N := \opn{Ker}(\phi^*)$
is a free left $\k \bra{z_1, z_2}$-module with basis
$\{ c s \mid s \in S \}$, where $c := z_1 z_2 - z_2 z_2$ and 
$S$ is the set of monomials.  
It is known that the category of finitely presented 
$\k \bra{z_1, z_2}$-modules (sometimes called coherent modules) is 
an abelian subcategory of 
$\cat{Mod} \k \bra{\bsym{z}}$; see \cite[Appendix, Theorem A.9]{Co}. 
Therefore from the exact sequence 
\[ 0 \to N \to \k \bra{z_1, z_2} \xar{\phi^*} \k[z_1, z_2]
\to 0 \]
we see that $\k[z_1, z_2]$ is not a finitely presented 
$\k \bra{z_1, z_2}$-module. It follows that there does 
not exists any exact sequence
\[ \k \bra{z_1, z_2}^{r_1} \to \k \bra{z_1, z_2}^{r_0}
\to \k [z_1, z_2] \to 0 \]
with $r_0, r_1 \in \mbb{N}$.
\end{exa}

\end{document}